\documentstyle[12pt]{article}
    \title{{\bf  Intertwining operator algebras and 
vertex tensor categories for affine Lie algebras}}
    \author{Yi-Zhi Huang
and James Lepowsky}
    \date{}
    \begin{document}
    \bibliographystyle{alpha}
    \maketitle

    \input amssym.def
    \input amssym
    \newtheorem{rema}{Remark}[section]
    \newtheorem{propo}[rema]{Proposition}
    \newtheorem{theo}[rema]{Theorem}
   \newtheorem{defi}[rema]{Definition}
    \newtheorem{lemma}[rema]{Lemma}
    \newtheorem{corol}[rema]{Corollary}
     \newtheorem{exam}[rema]{Example}
\newcommand{\binom}[2]{{{#1}\choose {#2}}}
	\newcommand{\nno}{\nonumber}
	\newcommand{\lbar}{\bigg\vert}
	\newcommand{\p}{\partial}
	\newcommand{\dps}{\displaystyle}
	\newcommand{\bra}{\langle}
	\newcommand{\ket}{\rangle}
 \newcommand{\res}{\mbox{\rm Res}}
\renewcommand{\hom}{\mbox{\rm Hom}}
 \newcommand{\pf}{{\it Proof}.\hspace{2ex}}
 \newcommand{\epf}{\hspace{2em}$\Box$}
 \newcommand{\epfv}{\hspace{1em}$\Box$\vspace{1em}}
\newcommand{\nord}{\mbox{\scriptsize ${\circ\atop\circ}$}}
\newcommand{\wt}{\mbox{\rm wt}\ }

 \makeatletter
\newlength{\@pxlwd} \newlength{\@rulewd} \newlength{\@pxlht}
\catcode`.=\active \catcode`B=\active \catcode`:=\active \catcode`|=\active
\def\sprite#1(#2,#3)[#4,#5]{
   \edef\@sprbox{\expandafter\@cdr\string#1\@nil @box}
   \expandafter\newsavebox\csname\@sprbox\endcsname
   \edef#1{\expandafter\usebox\csname\@sprbox\endcsname}
   \expandafter\setbox\csname\@sprbox\endcsname =\hbox\bgroup
   \vbox\bgroup
  \catcode`.=\active\catcode`B=\active\catcode`:=\active\catcode`|=\active
      \@pxlwd=#4 \divide\@pxlwd by #3 \@rulewd=\@pxlwd
      \@pxlht=#5 \divide\@pxlht by #2
      \def .{\hskip \@pxlwd \ignorespaces}
      \def B{\@ifnextchar B{\advance\@rulewd by \@pxlwd}{\vrule
         height \@pxlht width \@rulewd depth 0 pt \@rulewd=\@pxlwd}}
      \def :{\hbox\bgroup\vrule height \@pxlht width 0pt depth
0pt\ignorespaces}
      \def |{\vrule height \@pxlht width 0pt depth 0pt\egroup
         \prevdepth= -1000 pt}
   }
\def\endsprite{\egroup\egroup}
\catcode`.=12 \catcode`B=11 \catcode`:=12 \catcode`|=12\relax
\makeatother

\def\hboxtr{\FormOfHboxtr} 
\sprite{\FormOfHboxtr}(25,25)[0.5 em, 1.2 ex] 

:BBBBBBBBBBBBBBBBBBBBBBBBB |
:BB......................B |
:B.B.....................B |
:B..B....................B |
:B...B...................B |
:B....B..................B |
:B.....B.................B |
:B......B................B |
:B.......B...............B |
:B........B..............B |
:B.........B.............B |
:B..........B............B |
:B...........B...........B |
:B............B..........B |
:B.............B.........B |
:B..............B........B |
:B...............B.......B |
:B................B......B |
:B.................B.....B |
:B..................B....B |
:B...................B...B |
:B....................B..B |
:B.....................B.B |
:B......................BB |
:BBBBBBBBBBBBBBBBBBBBBBBBB |

\endsprite

\def\shboxtr{\FormOfShboxtr} 
\sprite{\FormOfShboxtr}(25,25)[0.3 em, 0.72 ex] 

:BBBBBBBBBBBBBBBBBBBBBBBBB |
:BB......................B |
:B.B.....................B |
:B..B....................B |
:B...B...................B |
:B....B..................B |
:B.....B.................B |
:B......B................B |
:B.......B...............B |
:B........B..............B |
:B.........B.............B |
:B..........B............B |
:B...........B...........B |
:B............B..........B |
:B.............B.........B |
:B..............B........B |
:B...............B.......B |
:B................B......B |
:B.................B.....B |
:B..................B....B |
:B...................B...B |
:B....................B..B |
:B.....................B.B |
:B......................BB |
:BBBBBBBBBBBBBBBBBBBBBBBBB |

\endsprite

\renewcommand{\theequation}{\thesection.\arabic{equation}}
\renewcommand{\therema}{\thesection.\arabic{rema}}
\setcounter{equation}{0}
\setcounter{rema}{0}
\setcounter{section}{-1}

\section{Introduction}

The category of finite direct sums of standard (integrable highest
weight) modules of a fixed positive integral level $k$ for an affine
Lie algebra $\hat{\frak g}$ is particularly important {}from the
viewpoint of conformal field theory and related mathematics. Here we
call this the category {\it generated by} the standard
$\hat{\frak g}$-modules of level $k$.  A central theme is a braided
tensor category structure (in the sense of Joyal and Street \cite{JS})
on this category, a structure explicitly discovered by Moore and
Seiberg \cite{MS} in their important study of conformal field
theories. In \cite{MS}, Moore and Seiberg constructed this structure
based on the assumption that there exists a suitable operator product
expansion for chiral vertex operators; this is essentially equivalent
to assuming the associativity of intertwining operators, in the
language of vertex operator algebra theory. Actually, Belavin,
Polyakov and Zamolodchikov \cite{BPZ} had already formalized the
relation between the operator product expansion and representation
theory in the context of conformal field theory, especially for the
Virasoro algebra, and Knizhnik and Zamolodchikov \cite{KZ} had
established the relation between conformal field theory and the
representation theory of affine Lie algebras.

In \cite{KL1}--\cite{KL5}, Kazhdan and Lusztig achieved a breakthrough
by indeed constructing a natural braided tensor category structure,
with the additional property of rigidity, on a certain category of
$\hat{\frak g}$-modules of level $k$, when $k$ is sufficiently
negative, or more generally, when $k$ is in a certain large subset of
${\Bbb C}$ excluding the positive integers, and, particularly, proving
that this braided tensor category is equivalent to a tensor category
of modules for a quantum group constructed {}from the same
finite-dimensional Lie algebra. The method used by Kazhdan and
Lusztig, especially in their construction of the associativity
isomorphisms, is algebro-geometric and is closely related to the
algebro-geometric formulation and study of conformal-field-theoretic
structures in the influential works of Tsuchiya-Ueno-Yamada
\cite{TUY}, Drinfeld \cite{Dr} and Beilinson-Feigin-Mazur \cite{BFM}.
(The work \cite{BFM} discusses the case of the minimal models for the
Virasoro algebra, and Beilinson informs us that a similar argument
works for the case of the category generated by the standard
$\hat{\frak g}$-modules.)

In the important work \cite{F} and \cite{F2}, Finkelberg proved that the
braided tensor category at positive level is tensor equivalent to a
certain ``subquotient'' of the braided tensor category constructed by
Kazhdan and Lusztig at negative level and is thus equivalent to a
certain ``subquotient'' braided tensor category of quantum group
representations.  (Cf. also \cite{Va}.) We have been informed by
Finkelberg that the arguments in his paper \cite{F2} can in fact be
reinterpreted to actually give a proof of the coherence relations for
the braided tensor category structure on the category generated by the
standard $\hat{\frak g}$-modules of positive integral level $k$; this
reasoning uses Kazhdan-Lusztig's result mentioned above constructing
rigid braided tensor category structure at negative level, and it also uses
Deligne's reformulation \cite{De} of the notion of balanced braided
tensor category.  (Finkelberg also informs us that his arguments do
not establish the coherence in a few cases, namely, those of ${\frak g}=E_6,
E_7, E_8$, $k=1$ or ${\frak g}=E_8$, $k=2$, because of the
corresponding restrictions at negative level in
\cite{KL1}--\cite{KL5}.)

In this paper, we prove the associativity of intertwining operators
(the existence of operator product expansion for chiral vertex
operators) for the vertex operator algebras associated to affine Lie
algebras at positive level, and we thereby construct directly the
braided tensor category structure on the category generated by the
standard $\hat{\frak g}$-modules of positive integral level, as an
application of the general tensor product theory for representations
of a vertex operator algebra developed in our papers
\cite{HL1}--\cite{HL6} and \cite{H1}--\cite{H2}.  We also prove more
general results, and in fact the methods in our general theory are
very different from those in the works \cite{MS},
\cite{KL1}--\cite{KL5}, \cite{TUY}, \cite{Dr}, \cite{BFM} and
\cite{F}, \cite{F2} mentioned above. We hope that the present work
helps illustrate the viewpoint that vertex operator algebra theory is
the appropriate mathematical foundation for understanding the range of
conformal-field-theoretic structures, including the associativity of
intertwining operators and the braided tensor category structure.

What the present paper specifically does is to use the works of
Knizhnik-Zamolodchikov \cite{KZ}, Tsuchiya-Kanie \cite{TK},
Frenkel-Huang-Lepowsky \cite{FHL}, Frenkel-Zhu \cite{FZ},
Dong-Li-Mason \cite{DLM}, Dong-Mason \cite{DM}, Dong-Mason-Zhu
\cite{DMZ} and Li \cite{L} \cite{L2} to verify the technical
conditions needed to apply our tensor product theory in the particular
case of the conformal field theories associated with affine Lie
algebras (or Wess-Zumino-Novikov-Witten models) and related models in
conformal field theory; the main part of our construction is already
contained in \cite{HL1}--\cite{HL6} and especially
\cite{H1}--\cite{H2}.  We would like to emphasize that the fundamental
ideas in the works Knizhnik-Zamolodchikov \cite{KZ}, Tsuchiya-Kanie
\cite{TK} and Frenkel-Zhu \cite{FZ} play essential roles in the
present paper, as does the work \cite{FHL}.

Our general theory in \cite{HL1}--\cite{HL6} and \cite{H1}-\cite{H2} (in
the generality of suitable vertex operator algebras) was initially
motivated by Kazhdan's and Lusztig's major series of papers
\cite{KL1}--\cite{KL5}.  The category studied by Kazhdan and Lusztig
is significantly different {}from the category generated by the
standard $\hat{\frak g}$-modules of level $k$.  It is larger; it
consists of the modules of level $k$ of finite length and whose
composition factors are irreducible highest-weight modules
corresponding to weights that are dominant integral in the direction
of $\frak g$.  For $k$ in the range studied by Kazhdan and Lusztig,
there are no nonzero standard $\hat{\frak g}$-modules of level $k$,
and for $k\in {\Bbb Z}_{+}$, there is no braided tensor category
structure on the category considered by Kazhdan and Lusztig (as is
pointed out in their work).  Correspondingly, our methods are very
different {}from Kazhdan's and Lusztig's, with an important exception:
The {\it definition} of the tensor product operation (for general
vertex operator algebras) in our work (see \cite{HL1}--\cite{HL4}) is
analogous to and was inspired by that in \cite{KL1}--\cite{KL5}.
However, our approach is based on a certain determination of the
tensor product module (see \cite{HL1}--\cite{HL5}).  Moreover, the
{\it construction of the associativity isomorphisms} (see \cite{H1})
uses a completely different method {}from the one in
\cite{KL1}--\cite{KL5}; a large amount of the information needed in
the proof of the necessary coherence properties (i.e., the pentagon
and hexagon properties) is already encoded in the associativity
isomorphisms.

The associativity of intertwining operators for the vertex operator
algebra associated to an affine Lie algebra at positive level proved
in this paper implies that the direct sum of all inequivalent
irreducible modules for the vertex operator algebra forms an
intertwining operator algebra, in the sense of \cite{H2.5}, \cite{H4}
and \cite{H6}.  We are also able in fact to obtain what we have called a
``vertex tensor category'' structure (see \cite{HL4} and \cite{HL6}),
which is
much richer than a braided tensor category structure; it is defined in terms
of certain moduli spaces of genus-zero Riemann surfaces with punctures
and local coordinates, equipped with the sewing operation (see \cite{H3}).
Vertex tensor category structure automatically yields braided tensor
category structure by the process of ignoring the conformal-geometric
information and retaining only the topological information.  We
further obtain the corresponding results for certain vertex operator
algebras slightly more general than those associated with the WZNW
models.  As we said above, we do this in this paper by verifying the
technical conditions needed to apply our general theory.

It was proved in \cite{H1} that the existence of the associativity
isomorphisms is equivalent to the associativity of intertwining
operators, in the language of vertex operator algebra theory.  In
their fundamental paper \cite{TK}, Tsuchiya and Kanie studied (in the
special case ${\frak g}={\frak s}{\frak l}(2, {\Bbb C})$) correlation
functions constructed {}from products of several intertwining
operators, and they used these correlation functions to construct and
study braid group representations.  Their main tool was the classic
Knizhnik-Zamolodchikov differential equations \cite{KZ}.  In the
present paper, we use this same method to verify a basic technical
condition, namely, the convergence and extension property, a key condition
needed for the application of the main associativity theorem in
\cite{H1}.  However, the definition of intertwining operator (or
chiral vertex operator) in \cite{TK} is different {}from the one in
\cite{MS} or the one in \cite{FHL}. The definition in \cite{TK} (as in
many other works) is restricted to primary fields (i.e., intertwining
operators in the sense of \cite{FHL} associated with the lowest
conformal-weight vectors of a module), and the iterated action of
$\hat{\frak g}$ on these operators is also fundamentally used.  This
restricted definition of intertwining operator is intimately related
to the prevalent notion of correlation function based on certain Lie
algebra coinvariants (as in \cite{KL1}--\cite{KL5}, for example).  Our
notion of correlation function is instead based on products of
intertwining operators in the sense of \cite{FHL}.

The ``nuclear democracy theorem'' in \cite{TK} can in fact be
reinterpreted as asserting the equivalence of these two different
notions of intertwining operator in the case ${\frak g}={\frak
s}{\frak l}(2, {\Bbb C})$, but iterates of intertwining operators
themselves, essential to the formulation and proof of the
associativity of intertwining operators, were not studied in
\cite{TK}.  The proof of the results in \cite{TK} necessary for the
construction of the braid group representations used the KZ equations
as well as certain algebraic constraints based on null-vector
conditions.

The notion of intertwining operator in \cite{FHL} (which we use in the
present paper) is the notion that is natural {}from the general
perspective of vertex operator algebra theory. (The starting point is
the Jacobi identity for vertex operator algebras \cite{FLM}; cf.
\cite{B}.)  In fact, the theory in \cite{HL1}--\cite{HL6} and
\cite{H1}--\cite{H2}, based on this notion, automatically incorporates
the corresponding algebraic constraints and the subtlety of ``nuclear
democracy'' at every stage, and the theory 
works in considerable generality. This
notion has led us to consider iterates of intertwining operators, to
formulate the associativity of intertwining operators (in terms of
such iterates) and to formulate a ``compatibility condition,'' which
was used to construct our tensor products in \cite{HL2} and
\cite{HL3}.  On the other hand, it was shown by Li \cite{L} \cite{L3}
that the ``nuclear democracy theorem'' generalizes to any ${\frak g}$,
and thus for any ${\frak g}$, the notions of intertwining operator in
the sense of \cite{TK} and in the sense of \cite{FHL} are indeed
equivalent.

In \cite{H1}, the associativity of intertwining operators was proved
using our theory of tensor products of modules for a vertex operator
algebra developed in \cite{HL2}, \cite{HL3} and \cite{HL5} when the
vertex operator algebra is rational in the sense of \cite{HL2} and
satisfies certain additional technical conditions. Combined with the
results of Moore and Seiberg \cite{MS}, the result of \cite{H1} served
to construct a natural braided tensor category structure on the
category of modules for such a vertex operator algebra.  
We also showed in \cite{H2} and \cite{HL6} that when the vertex operator
algebra contains a rational vertex operator subalgebra (but may itself
not be rational) and satisfies the additional technical conditions
mentioned above, we still get a braided tensor category structure such
that the tensor product bifunctor is the bifunctor $\boxtimes_{P(1)}$
constructed in \cite{HL5}, and in fact we get a vertex tensor category
structure.  In \cite{H2}, it was verified that for a vertex operator
algebra containing a vertex operator subalgebra isomorphic to a tensor
product of minimal Virasoro vertex operator algebras, the technical
conditions for the applicability of the tensor product theory are
satisfied. Thus the category of modules for such a vertex operator
algebra has a natural vertex tensor category structure and in
particular, a natural braided tensor category structure.

In the present paper, we work in the following generality, somewhat
greater than the generality of vertex operator algebras associated
with $\hat{\frak g}$ and positive integral level $k$: We verify that
for a vertex operator algebra containing a vertex operator subalgebra
isomorphic to a tensor product of vertex operator algebras associated
with WZNW models, the technical conditions for the applicability of
our tensor product theory are satisfied. (Recall {}from \cite{FHL} the
notions of tensor product and subalgebra for vertex operator algebras;
in particular, a subalgebra is required to have the same Virasoro
element as the large algebra.)  We thereby complete the construction
of the desired vertex tensor category structure, and in particular,
braided tensor category structure, on the category of modules for such
a vertex operator algebra.  In the special case of a vertex operator
algebra associated with a WZNW model, we conclude in particular that
the category generated by the standard $\hat{\frak g}$-modules of
positive integral level $k$ is indeed a braided tensor category. See
Section 3 for the main results.

The main tool that we use in this paper to verify our conditions is
the classical system of Knizhnik-Zamolodchikov equations \cite{KZ}.
Since we work in the framework of the theory of vertex operator
algebras, we include an exposition in the language of vertex operator
algebras of the fundamental theorem in \cite{KZ} that products of
intertwining operators for a vertex operator algebra associated with
an affine Lie algebra satisfy the Knizhnik-Zamolodchikov equations.
Our treatment of this result is nothing but an adaptation of the
original argument in \cite{KZ}, except that we use the language of formal
variables, in the spirit of the present theory; cf.  the treatment of
the Knizhnik-Zamolodchikov equations in \cite{TK}.  As in \cite{KZ}
and \cite{TK}, the theorem that products of intertwining operators
satisfy the Knizhnik-Zamolodchikov equations must be understood as an
assertion about formal series and not about functions of complex
variables, since one does not yet know at this stage that these
products are convergent; in fact, the Knizhnik-Zamolodchikov equations
are used to {\it prove} this convergence of products of intertwining
operators.

We assume that the reader is familiar with the statements of the main
results of \cite{HL1}--\cite{HL6} and \cite{H1}-\cite{H2}. We use the
basics of formal calculus and vertex operator algebra theory as
presented in \cite{FHL}, for example.  We shall assume in this paper
that in the definition of module for a vertex operator algebra, the
grading is by $\Bbb C$, not $\Bbb Q$ (see Remark 4.1.2 in \cite{FHL}).
Also, recall that a module by definition has two grading-restriction
properties---finite-dimensionality of the weight spaces and lower
truncation of the grading, and in this paper, when we assert that a
structure is a module, these restrictions are often the main issue.

In Section 1, we briefly review vertex operator algebras associated
with affine Lie algebras and their representations.  The results in
this section are due to Frenkel and Zhu \cite{FZ} (see also \cite{DL}
and \cite{L2}).  In Section 2, we present our exposition that products
of intertwining operators for a vertex operator algebra associated
with an affine Lie algebra satisfy the Knizhnik-Zamolodchikov
equations \cite{KZ}.  In Section 3, we use the Knizhnik-Zamolodchikov
equations and results obtained in \cite{HL1}--\cite{HL6} and
\cite{H1}-\cite{H2}, combined with results in \cite{DLM}, \cite{DM},
\cite{DMZ}, \cite{FHL} and \cite{L}, to prove our main results.

\paragraph{Acknowledgments} We are very grateful to Sasha Beilinson,
Pierre Deligne, Michael Finkelberg and Sasha Kirillov for discussions
concerning braided tensor category structure on the category generated
by the standard modules of positive integral level for an affine Lie
algebra.  We would like to thank Haisheng Li for discussions on
rationality of vertex operator algebras.  We would also like to thank
the referees for their questions and comments, which we have addressed
in the Introduction.  Y.-Z.~H. is supported in part by NSF grant
DMS-9596101 and DMS-9622961 and J.~L. by NSF grant DMS-9401851.

\renewcommand{\theequation}{\thesection.\arabic{equation}}
\renewcommand{\therema}{\thesection.\arabic{rema}}
\setcounter{equation}{0}
\setcounter{rema}{0}

\section{Vertex operator algebras associated with affine Lie algebras}

Let ${\frak g}$ be a  Lie algebra over ${\Bbb
C}$ equipped with an invariant symmetric bilinear form $(\cdot, \cdot)$. 
The {\it affine Lie algebra} $\hat{\frak g}$ associated with ${\frak g}$
and $(\cdot, \cdot)$ 
is the vector space ${\frak g}\otimes [t, t^{-1}]\oplus {\Bbb C}{\bf
k}$ equipped with the bracket operation defined by
\begin{eqnarray*}
[a\otimes t^{m}, b\otimes t^{n}]&=&[a, b]\otimes t^{m+n}+(a,
b)m\delta_{m+n, 0}{\bf k},\nno\\
{[a\otimes t^{m}, {\bf k}]}&=&0,
\end{eqnarray*}
for $a, b\in {\frak g}$ and $m, n\in {\Bbb Z}$. It is ${\Bbb Z}$-graded
in a natural way. Consider the subalgebras
$$\hat{\frak g}_{\pm}={\frak g}\otimes t^{\pm 1}{\Bbb C}[t^{\pm 1}]
\quad\mbox{\rm and}\quad
{\frak g}={\frak g}\otimes {\Bbb C}  
\quad\mbox{\rm of}\quad
\hat{\frak g} $$
and the vector
space decomposition
$$\hat{\frak g}=\hat{\frak g}_{-}\oplus {\frak
g}\oplus {\Bbb C}{\bf k}\oplus \hat{\frak g}_{+}.$$

Let $M$ be a ${\frak g}$-module, viewed as homogeneously graded of 
fixed degree $l\in {\Bbb C}$,   and let $k\in {\Bbb C}$. 
Let  $\hat{\frak g}_{+}$ act  trivially on $M$ and
${\bf k}$ as the scalar multiplication operator $k$. Then
$M$ becomes a ${\frak g}\oplus {\Bbb C}{\bf
k}\oplus \hat{\frak g}_{+}$-module, and we have the ${\Bbb C}$-graded
induced $\hat{\frak g}$-module
$$\hat{M}_{k}=U(\hat{\frak g})\otimes_{U(
{\frak g}\oplus {\Bbb C}{\bf
k}\oplus \hat{\frak g}_{+})} M,$$
which contains a canonical copy of $M$, of degree $l$.

Now we assume that ${\frak g}$ is finite-dimensional and simple. Let
 ${\frak h}$ be a Cartan subalgebra, $\Delta$ the root system
of $({\frak g}, {\frak h})$, $\Delta_{+}\subset \Delta$ 
a fixed set of positive
roots, $\theta$ the highest root, $h\check{~}$ the 
dual Coxeter number of ${\frak g}$, and $(\cdot, \cdot):
{\frak g}\times {\frak g}\to {\Bbb C}$ the Cartan-Killing form,
normalized by the condition $(\theta, \theta)=2$, where $(\cdot, \cdot)$
is transported canonically to ${\frak h}^{*}\times {\frak h}^{*}$. 
For a fixed  $\lambda\in {\frak h}^{*}$, let $L(\lambda)$ be the
irreducible highest-weight ${\frak g}$-module with highest weight
$\lambda$, assigned a homogeneous degree in ${\Bbb C}$. 
We shall use the notation $M(k, \lambda)$ to denote the graded $\hat{\frak
g}$-module $\widehat{L(\lambda)}_{k}$ (the degree of $L(\lambda)$
being suppressed in the notation). Let $J(k, \lambda)$ be the maximal
proper graded submodule of $M(k, \lambda)$ and $L(k, \lambda)=M(k,
\lambda)/J(k, \lambda)$. Then $L(k, \lambda)$ is the unique
irreducible graded 
$\hat{\frak g}$-module such that ${\bf k}$ acts as $k$ and
the space of all elements annihilated by $\hat{\frak g}_{+}$ is
isomorphic to the ${\frak g}$-module $L(\lambda)$. Later we shall 
often use the same  notations to denote elements
 of $M(k, \lambda)$ and of $L(k, \lambda)$.

In the special case $\lambda=0$, $L(0)$ is the one-dimensional trivial
${\frak g}$-module (equipped with a homogeneous degree in ${\Bbb C}$)
 and can be identified with ${\Bbb C}$. Consequently,
as a vector space and $\hat{\frak
g}_{-}$-module, $M(k, 0)$ is naturally isomorphic to the universal
enveloping algebra $U(\hat{\frak g}_{-})$. We define a vertex operator
map 
\begin{eqnarray*}
Y(\cdot, x): M(k, 0)&\to& (\mbox{\rm End}\;M(k, 0))[[x, x^{-1}]]\nno\\
v&\mapsto &Y(v, x)=\sum_{n\in {\Bbb Z}}v_{n}x^{-n-1}
\end{eqnarray*}
as follows: Identifying $M(k, 0)$ with $U(\hat{\frak g}_{-})$, we note
that $M(k, 0)$ is spanned by the elements of the form
$a_{1}(-n_{1})\cdots a_{m}(-n_{m})1$, where $a_{1}, \dots, a_{m}\in
{\frak g}$ and $n_{1}, \dots, n_{m}\in {\Bbb Z}_{+}$, with $a(-n)$ denoting
 the representation image of
$a\otimes t^{-n}$ for $a\in {\frak g}$ and $n\in {\Bbb Z}$ 
(we shall use similar notation for other $\hat{\frak
g}$-modules below). We use recursion on $m$ to define the vertex
operator map. For $1\in M(k, 0)$, we define $Y(1, x)$ to be the
identity operator on $M(k, 0)$ and for $a\in {\frak g}$ we define 
$$Y(a(-1)1, x)=\sum_{n\in {\Bbb Z}}a(n)x^{-n-1}.$$
Assume that the  map has been
 defined for  $a_{1}(-n_{1})\cdots
a_{m}(-n_{m})1$. Vertex operators for elements
of the form $a_{0}(-n_{0})a_{1}(-n_{1})\cdots a_{m}(-n_{m})1$ are defined 
using the residue in $x_{1}$
of the Jacobi identity for 
vertex operator algebras (cf. (\ref{3.0}) below) as follows:
\begin{eqnarray}\label{2.0}
\lefteqn{Y(a_{0}(-n_{0})a_{1}(-n_{1})\cdots a_{m}(-n_{m})1, x)=}\nno\\
&&=\res_{x_{1}}(x_{1}-x)^{-n_{0}}Y(a_{0}(-1)1, x_{1})
Y(a_{1}(-n_{1})\cdots a_{m}(-n_{m})1, x)\nno\\
&&\quad -\res_{x_{1}}(-x+x_{1})^{-n_{0}}
Y(a_{1}(-n_{1})\cdots a_{m}(-n_{m})1, x)Y(a_{0}(-1)1, x_{1}).\nno\\
&&
\end{eqnarray}
(Here and below, binomial expressions are understood to be 
expanded in nonnegative powers of the second variable.)
The vacuum vector for  $M(k, 0)$ is ${\bf 1}=1$. In the case that $k\ne -
h\check{~}$,
$M(k, 0)$ also has a Virasoro element
\begin{equation}\label{omega}
\omega=\frac{1}{2(k+h\check{~})}\sum_{i=1}^{\dim {\frak g}}
g^{i}(-1)^{2}{\bf 1},
\end{equation}
where $\{g^{i}\}_{i=1, \dots, \dim {\frak g}}$ is an arbitrary orthonormal
basis of ${\frak g}$ with respect to the  form $(\cdot, \cdot)$.
With respect to the associated grading, ${\bf 1}$ has weight $0$; the
(conformal) weight of a homogeneous element is the negative of its degree.
In \cite{FZ}, the following result is proved (among other things):

\begin{theo}\label{2-1}
If $k\ne -h\check{~}$, the quadruple 
$(M(k, 0), Y, {\bf 1}, \omega)$ defined above is a vertex operator algebra;
in particular, $Y(\cdot, x)$ is well defined.
\end{theo}

 Since
 $J(k, 0)$ is a $\hat{\frak g}$-submodule of $M(k, 0)$, 
the vertex operator map for $M(k, 0)$ induces a vertex 
operator map for $L(k, 0)$ which we shall still denote 
$Y$. We continue to denote the $J(k, 0)$-cosets of ${\bf 1}$ and $\omega$ in 
$M(k, 0)$ by ${\bf 1}$ and $\omega$.
The following result is an immediate consequence
of Theorem \ref{2-1}:

\begin{corol}[\cite{FZ}]\label{2-2}
If
$k\ne -h\check{~}$, the quadruple 
$(L(k, 0), Y, {\bf 1}, \omega)$ is a (non\-zero) vertex operator algebra.
\end{corol}

We now discuss modules for these vertex operator algebras. 
For any $k\in {\Bbb C}$ and $\lambda\in {\frak h}^{*}$, we define 
a vertex operator map 
$$Y(\cdot, x): M(k, 0)\to (\mbox{\rm End}\ M(k, \lambda))[[x, x^{-1}]]$$
using recursion, just as in (\ref{2.0}). As in Theorem \ref{2-1}, this
is well defined.
The induced map
$$Y(\cdot, x): L(k, 0)\to (\mbox{\rm End}\ L(k, \lambda))[[x,
x^{-1}]]$$
is well defined under the conditions given in 
the following result, also proved in \cite{FZ}; see also \cite{L2}:

\begin{theo}\label{fz}
For any $k\in {\Bbb C}$ such that $k\ne -h\check{~}$ and any
$\lambda\in {\frak h}^{*}$ with $\lambda$ dominant integral
(i.e., with
$\dim L(\lambda)<\infty$),
the pair $(M(k, \lambda), Y)$ is a  module for
the vertex operator 
algebra
$M(k, 0)$. In case  $k=0, 1, 2,\dots$,
$L(k, 0)$ is a rational vertex operator algebra and 
$$\{L(k, \lambda)\;|\;\lambda\in {\frak h}^{*} \;\;\mbox{\rm is 
dominant integral such that}\;\; (\lambda,\theta)\le k\}$$
is the set of all irreducible $L(k, 0)$-modules up to equivalence.
\end{theo}

\begin{rema}
{\rm In the present paper, {\it rationality} for a vertex operator algebra
refers to the definition in \cite{HL2}: the number of
inequivalent irreducible modules is finite, every module is
completely reducible, and the fusion rules (the dimensions of the
spaces of intertwining operators among triples of modules) are finite.
The notion of rationality used in \cite{FZ} was a different one,
but the rationality of $L(k,0)$ in our sense in fact holds.}
\end{rema}
 
\begin{rema}
{\rm When $k\in {\Bbb N}$,
the vertex operator algebras $L(k, 0)$ and their irreducible modules can 
also be constructed explicitly using Fock spaces and tensor products. 
See \cite{B}, \cite{FLM} and Chapter 13 of \cite{DL}.}
\end{rema}

\renewcommand{\theequation}{\thesection.\arabic{equation}}
\renewcommand{\therema}{\thesection.\arabic{rema}}
\setcounter{equation}{0}
\setcounter{rema}{0}

\section{Products of intertwining operators and the KZ equations}

We give our exposition in this section of the following:  When $k\in
{\Bbb N}$, products of intertwining operators for the vertex operator
algebra $L(k, 0)$ satisfy a system of differential equations with
regular singular points called the {\it Knizhnik-Zamolodchikov
equations} (or simply the {\it KZ equations}), which were first
derived in \cite{KZ}.  The reason for our supplying this exposition
was explained in the Introduction.  The discussions and results in
this section also hold for the vertex operator algebras $M(k, 0)$ for
any $k\in {\Bbb C}$ such that $k\ne -h\check{~}$.

In this  and the next section, we shall need a simple general
principle for vertex operators and intertwining operators
(see (\ref{3.1}) below; cf. \cite{DL}, formulas (13.24)--(13.26))). 
We shall use the basics of formal calculus as expressed in Section 2.1
of \cite{FHL}, in particular, the binomial expansion convention. 

 Let $V$ be a
vertex operator algebra, let $W_{1}, W_{2}, W_{3}$ be 
$V$-modules and let ${\cal
Y}$ be an intertwining operator of type ${W_{3}\choose W_{1}W_{2}}$. Then
{}from the Jacobi identity defining ${\cal Y}$ (see
\cite{FHL}, formula (5.4.4)), we obtain, by taking $\res_{x_{1}}$,
\begin{eqnarray}\label{3.0}
{\cal Y}(Y(u, x_{0})w, x_{2})&=&
\left(\res_{x_{1}}x_{0}^{-1}\delta\left(\frac{x_{1}-x_{2}}{x_{0}}\right)
Y(u, x_{1})\right){\cal Y}(w, x_{2})\nno\\
&& -{\cal Y}(w, x_{2})
\left(\res_{x_{1}}x_{0}^{-1}\delta\left(\frac{x_{2}-x_{1}}{-x_{0}}\right)
Y(u, x_{1})\right)
\end{eqnarray}
for $u\in V$ and $w\in W_{1}$ (cf. (\ref{2.0})).

We shall write
\begin{equation}\label{3.0.0.1}
Y(u, x)=Y^{+}(u, x)+Y^{-}(u, x)
\end{equation}
 where $Y^{+}(u, x)=\sum_{n<0}
u_{n}x^{-n-1}$ and $Y^{-}(u, x)=\sum_{n\ge 0}u_{n}x^{-n-1}$ 
are the regular and singular parts of $Y(u, x)$,
respectively. 
Equating the  parts of both sides of
(\ref{3.0}) which are constant in $x_{0}$, we obtain 
\begin{eqnarray}\label{3.1}
{\cal Y}(u_{-1}w, x_{2})
&=&\left(\res_{x_{1}}
(x_{1}-x_{2})^{-1}
Y(u, x_{1})\right){\cal Y}(w, x_{2})\nno\\
&&\quad +{\cal Y}(w, x_{2})
\left(\res_{x_{1}}
(x_{2}-x_{1})^{-1}
Y(u, x_{1})\right)\nno\\
&=&\left(\res_{x_{1}}
x_{1}^{-1}\delta\left(\frac{x_{2}}{x_{1}}\right)
Y^{+}(u, x_{1})\right){\cal Y}(w, x_{2})\nno\\
&&\quad +{\cal Y}(w, x_{2})
\left(\res_{x_{1}}x_{1}^{-1}\delta\left(\frac{x_{2}}{x_{1}}\right)
Y^{-}(u, x_{1})\right)\nno\\
&=&\left(\res_{x_{1}}
x_{1}^{-1}\delta\left(\frac{x_{2}}{x_{1}}\right)
Y^{+}(u, x_{2})\right){\cal Y}(w, x_{2})\nno\\
&&\quad +{\cal Y}(w, x_{2})
\left(\res_{x_{1}}x_{1}^{-1}\delta\left(\frac{x_{2}}{x_{1}}\right)
Y^{-}(u, x_{2})\right)\nno\\
&=&Y^{+}(u, x_{2}){\cal Y}(w, x_{2})
+{\cal Y}(w, x_{2})Y^{-}(u,
x_{2})\nno\\
&=&\nord Y(u, x_{2}){\cal Y}(w, x_{2})\nord,
\end{eqnarray}
where we define
the ``normal ordering'' $\nord \cdot\nord$  by
\begin{equation}\label{nord}
\nord u_{m}w_{n}\nord 
=\left\{ \begin{array}{ll}u_{m}w_{n}&\mbox{\rm if}\;\;
m<0\\w_{n}u_{m}&\mbox{\rm if}\;\; 
m\ge 0.\end{array}
\right.
\end{equation}
Note that this particular normal-ordering operation is 
not in general 
commutative; moreover, the expansion ${\cal Y}(w, x)=\sum_{n\in {\Bbb C}}
w_{n}x^{-n-1}$ involves nonintegral subscripts $n$ in general.

Now we take $V$ to be one of the vertex operator algebras $L(k, 0)$
for $k\in {\Bbb N}$ or $M(k, 0)$ for $k\in {\Bbb C}$, $k\ne
-h\check{~}$. Let $W$ be any $V$-module. 
 For any $g\in {\frak g}$, we denote $Y(g(-1){\bf 1}, x)$,
$Y^{+}(g(-1){\bf 1}, x)$ and $Y^{-}(g(-1){\bf 1}, x)$, acting on $W$,  by 
$g(x)$, $g^{+}(x)$ and $g^{-}(x)$, 
respectively.

{}From (\ref{omega}) and (\ref{3.1}) (note that 
$g^{i}(-1)=(g^{i}(-1){\bf 1})_{-1}$), we obtain 
\begin{eqnarray}
Y(\omega, x)&=&\frac{1}{2(k+h\check{~})}\sum_{i=1}^{\dim {\frak g}}
: g^{i}(x)^{2}:,\\
L(n)&=&\frac{1}{2(k+h\check{~})}\sum_{i=1}^{\dim {\frak g}}
\sum_{j\in {\Bbb Z}}:g^{i}(j)g^{i}(n-j):
\end{eqnarray}
for $n\in {\Bbb Z}$, where for any $g, g'\in {\frak g}$,
$$:g(r)g'(s):
=\left\{ \begin{array}{ll}g(r)g'(s)&\mbox{\rm if}\;\;
r< s\\
\frac{1}{2}(g(r)g'(s)+g'(s)g(r))&\mbox{\rm if}\;\;
r= s\\
g'(s)g(r)&\mbox{\rm if}\;\; r>s\end{array}
\right.$$
(cf. Proposition 13.4 and  (13.35) in \cite{DL}; in the present situation,
the two normal orderings coincide).
In particular,
\begin{eqnarray*}
L(-1)&=&\frac{1}{2(k+h\check{~})}\sum_{i=1}^{\dim {\frak g}}
\left(\sum_{j< 0}g^{i}(j)g^{i}(-j-1)
+\sum_{j\ge 0}g^{i}(-j-1)g^{i}(j)\right).
\end{eqnarray*}
If $w\in W$ is such that $g(n)w=0$
for $n>0$ and $g\in {\frak g}$, then
\begin{equation}\label{3.2}
L(-1)w=\frac{1}{k+h\check{~}}\sum_{i=1}^{\dim {\frak g}}
g^{i}(-1)g^{i}(0)w=
\frac{1}{k+h\check{~}}\sum_{i=1}^{\dim {\frak g}}g^{i}(-1)g^{i}w.
\end{equation}

Let $n\ge 1$, let $W_{0}, \dots, W_{n-1}$  be $V$-modules, and let
$W_{n}=V$. We shall consider a product of intertwining operators
of the form 
\begin{equation}\label{prod}
{\cal Y}_{1}(w_{(1)}, x_{1})\cdots {\cal Y}_{n}(w_{(n)}, x_{n})
\end{equation}
for $n\ge 1$.  Here $w_{(l)}\in L(\lambda_{l})$, where 
$\lambda_{l}\in {\frak h}^{*}$ is dominant integral. 
Recall Theorem \ref{fz}. In the case in which $V=L(k, 0)$ ($k\in {\Bbb
N}$), we assume that $(\lambda_{l},\theta)\le k$ and view
$L(\lambda_{l})$ as the lowest conformal-weight space of the $V$-module $L(k,
\lambda_{l})$, and in the case $V=M(k, 0)$ ($k\ne -h\check{~}$), we
analogously take $L(\lambda_{l})$ to be the lowest 
conformal-weight space of the
$V$-module $M(k, \lambda_{l})$. The intertwining operator 
${\cal Y}_{l}$ is
 of type 
$\binom{W_{l-1}}{L(k, \lambda_{l})\; W_{l}}$ in the former case and 
of type  $\binom{W_{l-1}}{M(k, \lambda_{l})\; W_{l}}$ in the latter case.

By the $L(-1)$-derivative property in the definition of intertwining
operator and (\ref{3.2}),
\begin{eqnarray}\label{3.3}
(k+h\check{~})\frac{d}{dx_{l}}{\cal Y}_{l}(w_{(l)}, x_{l})&=&(k+h\check{~})
{\cal Y}_{l}(L(-1)w_{(l)}, x_{l})\nno\\
&=&\sum_{i=1}^{\dim {\frak g}}\nord g^{i}(x_{l}){\cal Y}_{l}(g^{i}
w_{(l)}, x_{l})\nord.
\end{eqnarray}

For any $g\in {\frak g}$ and $w_{(p)}\in L(k, \lambda_{p})$ with $p\ne l$, 
\begin{eqnarray*}
\lefteqn{[g(x_{l}), {\cal Y}_{p}(w_{(p)}, x_{p})]}\nno\\
&&=\res_{x_{0}}x_{p}^{-1}\delta\left(\frac{x_{l}-x_{0}}{x_{p}}\right)
{\cal Y}_{p}(g(x_{0})w_{(p)}, x_{p})\nno\\
&&=\res_{x_{0}}x_{p}^{-1}\delta\left(\frac{x_{l}-x_{0}}{x_{p}}\right)
{\cal Y}_{p}(\sum_{n\le 0}g(n)x_{0}^{-n-1}w_{(p)}, x_{p})\nno\\
&&=x_{p}^{-1}\delta\left(\frac{x_{l}}{x_{p}}\right)
{\cal Y}_{p}(gw_{(p)}, x_{p}),
\end{eqnarray*}
and so
\begin{eqnarray}
[g^{-}(x_{l}),
{\cal Y}_{p}(w_{(p)}, x_{p})]
&=&(x_{l}-x_{p})^{-1}{\cal Y}_{p}(gw_{(p)}, x_{p}),\label{3.7.1}\\
{[g^{+}(x_{l}), {\cal Y}_{p}(w_{(p)}, x_{p})]}
&=&-(-x_{p}+x_{l})^{-1}{\cal Y}_{p}(gw_{(p)}, x_{p}),\label{3.7.2}
\end{eqnarray}
where we continue to use the binomial expansion 
convention.

We introduce the following notation: For  commuting 
formal variables $x_{1}, \dots, x_{n}$, let
$$:(x_{l}-x_{p})^{-1}:\;=\left\{\begin{array}{ll}(x_{l}-x_{p})^{-1}&\mbox{\rm
if} \;\;l<p\\
(-x_{p}+x_{l})^{-1}&\mbox{\rm
if} \;\;p<l.\end{array}\right.$$

Consider the contragredient $V$-module $W'_{0}$ (see \cite{FHL},
Sections 5.2 and 5.3) and let $w'_{(0)}\in W'_{0}$ be a lowest
conformal-weight vector.
By (\ref{3.3})--(\ref{3.7.2}) and since
\begin{eqnarray}
(g^{i})^{-}(x_{l}){\bf 1}&=&0,\nno\\
\langle w'_{(0)}, (g^{i})^{+}(x_{l})u\rangle&=&0\label{w'0}
\end{eqnarray}
for any $u\in V$,  we have (for fixed $l$)
\begin{eqnarray}\label{3.9}
\lefteqn{(k+h\check{~})\frac{d}{dx_{l}}\langle w'_{(0)}, {\cal Y}_{1}(w_{(1)}, 
x_{1})\cdots 
{\cal Y}_{n}(w_{(n)}, x_{n}){\bf 1}\rangle}\nno\\
&&=\sum_{p\ne l}:(x_{l}-x_{p})^{-1}:\sum_{i=1}^{\dim {\frak g}}
\langle w'_{(0)}, {\cal Y}_{1}(w_{(1)}, x_{1})\cdots 
{\cal Y}_{l}(g^{i}w_{(l)}, x_{l})\cdots\nno\\
&&\hspace{6em}\cdots
{\cal Y}_{p}(g^{i}w_{(p)}, x_{p})
\cdots
{\cal Y}_{n}(w_{(n)}, x_{n}){\bf 1}\rangle.
\end{eqnarray}

We define operators 
$\Omega_{lp}$, for $1\le l, p\le n$ with $l\ne p$, on 
$$(L(\lambda_{1})\otimes \cdots \otimes 
L(\lambda_{n}))^{*}$$
 as follows: 
$$(\Omega_{lp}f)(w_{(1)}\otimes \cdots \otimes  w_{(n)})
=\sum_{i=1}^{\dim {\frak g}}f(w_{(1)}\otimes \cdots \otimes 
g^{i}w_{(l)}\otimes \cdots\otimes 
g^{i}w_{(p)}\otimes \cdots \otimes w_{(n)}).$$
Note that $\Omega_{lp}=\Omega_{pl}$.
We extend $\Omega_{lp}$ naturally to an operator on the vector space 
$$(L(\lambda_{1})\otimes \cdots \otimes 
L(\lambda_{n}))^{*}\{x_{1}, \dots, x_{n}\}$$
 of formal series with arbitrary powers of the variables and with
coefficients in  $(L(\lambda_{1})\otimes \cdots \otimes 
L(\lambda_{n}))^{*}$.
In terms of these operators,  (\ref{3.9})
becomes:

\begin{theo}[Knizhnik-Zamolodchikov]\label{kzeq}
Let $V=L(k, 0)$ for $k\in {\Bbb
N}$ or $V=M(k, 0)$ for $k\in {\Bbb C}$, $k\ne -h\check{~}$. Let
$n\ge 1$, let $W_{0}, \dots, W_{n-1}$  be $V$-modules, and let
$W_{n}=V$. Let ${\cal Y}_{1}, \dots, {\cal Y}_{n}$ be intertwining
operators,
$w_{(l)}\in L(\lambda_{l})$ and $w'_{(0)}\in W'_{0}$ as described
above (see (\ref{prod}) and (\ref{w'0})), and define
$$\phi_{{\cal Y}_{1}, \dots, {\cal Y}_{n}}\in (L(\lambda_{1})\otimes
\cdots \otimes 
L(\lambda_{n}))^{*}\{x_{1}, \dots, x_{n}\}$$
 by 
$$\phi_{{\cal Y}_{1}, \dots, {\cal Y}_{n}}(w_{(1)}\otimes\cdots \otimes
w_{(n)})=\langle w'_{(0)},
{\cal Y}_{1}(w_{(1)}, x_{1})\cdots {\cal Y}_{n}(w_{(n)}, 
x_{n}){\bf 1}\rangle.$$
Then for $l=1, \dots, n$,
\begin{equation}\label{kzn}
(k+h\check{~})\frac{\p}{\p x_{l}}\phi_{{\cal Y}_{1}, \dots, {\cal Y}_{n}}=
\sum_{p\ne l}
:(x_{l}-x_{p})^{-1}:\Omega_{lp}\phi_{{\cal Y}_{1}, \dots, {\cal
Y}_{n}}.\;\;\;\Box
\end{equation}
\end{theo}

The equations (\ref{kzn}) are the {\it Knizhnik-Zamolodchikov
equations} (the {\it KZ equations}).  They are easily seen to be
consistent. Theorem \ref{kzeq} gives:

\begin{corol}
Let $k\in {\Bbb N}$; let 
$\lambda_{l}\in {\frak h}^{*}$, for $l=0, 1, 2, 3$, be dominant integral
weights satisfying $(\lambda_{l},\theta)\le k$; let $W$ be an $L(k, 0)$-module;
and let ${\cal Y}_{1}$ and ${\cal Y}_{2}$ be intertwining
operators 
of types 
$\binom{L(k, \lambda_{0})}{L(k, \lambda_{1})\; W}$ and 
$\binom{W}{L(k, \lambda_{2})\;L(k, \lambda_{3})}$, respectively. Define
$$\psi_{{\cal Y}_{1}, {\cal Y}_{2}}\in 
(L(\lambda_{0})^{*}\otimes L(\lambda_{1})\otimes 
L(\lambda_{2})\otimes L(\lambda_{3}))^{*}\{x_{1}, x_{2}\}$$ 
by
$$\psi_{{\cal Y}_{1}, {\cal Y}_{2}}(w'_{(0)}
\otimes w_{(1)}\otimes w_{(2)}\otimes w_{(3)})
=\langle w'_{(0)}, {\cal Y}_{1}(w_{(1)}, x_{1}){\cal Y}_{2}(w_{(2)}, x_{2})
w_{(3)}\rangle$$
for $w'_{(0)}\in L(\lambda_{0})^{*}\subset L(k, \lambda_{0})'$ and
$w_{(l)}\in L(\lambda_{l})\subset L(k, \lambda_{l})$, $l=1, 2, 3$.
Then 
\begin{eqnarray}
(k+h\check{~})\frac{\p}{\p x_{1}}\psi_{{\cal Y}_{1}, {\cal Y}_{2}}&=&
\left(\frac{\Omega_{13}}{x_{1}}
+\frac{\Omega_{12}}{x_{1}-x_{2}}\right)\psi_{{\cal Y}_{1}, 
{\cal Y}_{2}},\label{3.10}\\
(k+h\check{~})\frac{\p}{\p x_{2}}\psi_{{\cal Y}_{1}, {\cal Y}_{2}}&=&
\left(\frac{\Omega_{23}}{x_{2}}
+\frac{\Omega_{12}}{-x_{1}+x_{2}}\right)\psi_{{\cal Y}_{1}, {\cal
Y}_{2}}.\;\;\;\Box
\label{3.11}
\end{eqnarray}
\end{corol}

\renewcommand{\theequation}{\thesection.\arabic{equation}}
\renewcommand{\therema}{\thesection.\arabic{rema}}
\setcounter{equation}{0}
\setcounter{rema}{0}

\section{The category of 
 standard modules of a fixed level for an affine Lie algebra}

In this section we show that for the category of modules for a vertex
operator algebra containing a subalgebra isomorphic to a tensor
product of vertex operator algebras of the form $L(k, 0)$, $k\in {\Bbb
N}$, the intertwining operators among the modules have the
associativity property, the category has a natural structure of vertex
tensor category, and a number of other important results hold.  For
the tensor product theory for modules for a vertex operator algebra,
see \cite{HL1}--\cite{HL6} and \cite{H1}--\cite{H2}.

We need the following notion introduced in \cite{H1}:

\begin{defi}
{\rm Let $V$ be a vertex operator algebra. We say that products of
intertwining operators for $V$ satisfy the {\it convergence and
extension property} if for any intertwining operators ${\cal Y}_{1}$
and ${\cal Y}_{2}$ of types ${W_{0}}\choose {W_{1}W_{4}}$ and
${W_{4}}\choose {W_{2}W_{3}}$, respectively, there exists an integer
$N$ (depending only on ${\cal Y}_{1}$ and ${\cal Y}_{2}$)
such that the following condition holds:  For any
$w_{(1)}\in W_{1}$, $w_{(2)}\in W_{2}$, $w_{(3)}\in W_{3}$,
$w'_{(0)}\in W'_{0}$ with $w_{(1)}$ and $w_{(2)}$ homogeneous, there
exist $j\in {\Bbb N}$, $r_{i}, s_{i}\in {\Bbb R}$, for $i=1, \dots, j$,
and analytic functions $f_{i}(z)$ on $|z|<1$, for $i=1, \dots, j$,
satisfying
\begin{equation}\label{si}
\wt w_{(1)}+\wt w_{(2)}+s_{i}>N,\;\;\;i=1, \dots, j,
\end{equation}
such that
\begin{equation}\label{cor-fn}
\langle w'_{(0)}, {\cal Y}_{1}(w_{(1)}, x_{1})
{\cal Y}_{2}(w_{(2)}, x_{2})w_{(3)}\rangle
\lbar_{x_{i}^{r}=e^{r\log z_{i}}, \; i=1,2, \; r\in {\Bbb C}}
\end{equation}
is absolutely 
convergent when $|z_{1}|>|z_{2}|>0$ and can be analytically extended to  
the multivalued analytic function
\begin{equation}\label{phyper}
\sum_{i=1}^{j}z_{2}^{r_{i}}(z_{1}-z_{2})^{s_{i}}
f_{i}\left(\frac{z_{1}-z_{2}}{z_{2}}\right)
\end{equation}
when $|z_{2}|>|z_{1}-z_{2}|>0$;
here ``$\log$'' denotes the standard branch of the $\log$ function.}
\end{defi}

Note that if the associativity of intertwining operators holds for $V$
(see \cite{H1}-\cite{H2} and Theorem \ref{assoc} below), then products
of intertwining operators satisfy the convergence and extension
property, so that this condition is necessary for the associativity of
intertwining operators; to see that the expression in (\ref{1-1})
below has the form (\ref{phyper}), use the $L(0)$-conjugation formula
(5.4.22) in \cite{FHL}.

We shall also use the concepts of generalized module, as defined
in \cite{HL2}, and of weak module, as defined  in \cite{DLM}.
Let $V$ be a vertex operator algebra. A {\it generalized $V$-module}
 is a ${\Bbb C}$-graded vector space equipped with 
a vertex operator map satisfying all the axioms for a $V$-module except 
the two grading-restriction axioms. If there exists $N\in {\Bbb Z}$ 
such that 
the homogeneous subspace of weight $n$ of a generalized $V$-module
 is $0$ when the real part of 
$n$ is less than $N$, the generalized $V$-module is said to be {\it lower
truncated}. 

A {\it weak $V$-module} is a 
vector space equipped with 
a vertex operator map satisfying all the axioms for a $V$-module except 
those axioms involving a grading.
Dong, Li and Mason \cite{DLM} have proved a stronger result than
Theorem \ref{fz}---that for any $k\in {\Bbb N}$,
a weak $L(k, 0)$-module is
completely reducible and all the irreducible weak
$L(k, 0)$-modules are in fact $L(k, 0)$-modules, as listed in
Theorem \ref{fz}.  In particular, a finitely generated lower-truncated
generalized $L(k, 0)$-module is a module (see also \cite{FZ} and 
\cite{L2}, which use instead the notion of module in \cite{FZ}).

We also have that the tensor product vertex operator
algebra $L(k_{1}, 0)\otimes \cdots \otimes 
L(k_{m}, 0)$, for any $k_{1}, \dots, k_{m}\in {\Bbb N}$,
is rational, by the stronger result just stated, 
and by arguments in Section 4.7 of \cite{FHL}
and in Section 2 of \cite{DMZ}; see also 
\cite{L}.  Actually, for the complete reducibility of a module for the
tensor product vertex operator algebra, one can alternatively argue as
in \cite{FZ}, \cite{L2} and \cite{DMZ} using the notion of module in
\cite{FZ} rather than the notion of weak module.  (An argument
essentially the same as this complete reducibility argument appears in
fact in the proof of condition 1 of Theorem 3.2 below.)

The following result establishes the remaining conditions
for the applicability of the tensor product theory 
(see \cite{HL1}--\cite{HL6} and \cite{H1}--\cite{H2})
to the vertex operator algebra $L(k_{1}, 0)\otimes \cdots \otimes 
L(k_{m}, 0)$:

\begin{theo}\label{kim}
For any $m\in {\Bbb Z}_{+}$ and $k_{i}\in {\Bbb N}$, 
with $i=1, \dots, m$,  we have:

\begin{enumerate}

\item Every finitely-generated lower-truncated generalized module for
the vertex operator algebra 
$L(k_{1}, 0)\otimes \cdots \otimes 
L(k_{m}, 0)$ is a module.

\item Products of  intertwining operators for
$L(k_{1}, 0)\otimes \cdots \otimes 
L(k_{m}, 0)$ have the convergence and extension property.

\item For any modules $W_{j}$, with $j=1, \dots, 2n+1$,
for
$L(k_{1}, 0)\otimes \cdots \otimes 
L(k_{m}, 0)$; any  
intertwining operators ${\cal  Y}_{i}$, with $i=1, \dots, n$, 
of types ${W_{2i-1}}\choose {W_{2i}\; W_{2i+1}}$, respectively;
and any $w'_{(1)}\in W'_{1}$, $w_{(2i)}\in W_{2i}$, with $i=1, \dots, n$, 
and $w_{(2n+1)}\in W_{2n+1}$, 
\begin{equation}\label{cor-fn-n}
\langle w'_{(1)}, {\cal  Y}_{1}(w_{(2)}, x_{1})\cdots 
{\cal  Y}_{m}(w_{(2n)}, x_{n})w_{(2n+1)}\rangle\lbar_{x_{i}^{r}
=e^{r\log z_{i}}, \; 1\le i\le n, \; r\in {\Bbb C}}
\end{equation}
is absolutely convergent 
 for any $z_{1}, \dots, z_{n}\in {\Bbb C}$
satisfying $|z_{1}|>\cdots >|z_{n}|>0$; here the choice of $\log
z_{i}$ is arbitrary for each $i$.

\end{enumerate}
\end{theo}
\pf We first prove the theorem in the case $m=1$.  
We have already verified the first condition in this case.

To prove that products of intertwining operators for $L(k, 0)$ have
the convergence and extension property, we can assume that the $W_{l}$,
$l=0. \dots, 4$, are irreducible since $L(k, 0)$ is rational. Let
$\lambda_{l}\in {\frak h}^{*}$, $l=0, \dots, 3$, be dominant integral
weights satisfying $(\lambda_{l},\theta)\le k$, let $W$ be a $V$-module, and
let ${\cal Y}_{1}$ and ${\cal Y}_{2}$ be intertwining operators of types
$\binom{L(k, \lambda_{0})}{L(k, \lambda_{1})\; W}$ and $\binom{W}{L(k,
\lambda_{2})\; L(k, \lambda_{3})}$, respectively.  For $w'_{(0)}\in
L(\lambda_{0})^{*}\subset L(k, \lambda_{0})'$ and $w_{(i)}\in
L(\lambda_{i})\subset L(k, \lambda_{i})$, $i=1, 2, 3$, we have already
shown in the preceding section that (\ref{cor-fn}) satisfies the KZ
equations (\ref{3.10}) and (\ref{3.11}). 
Thus it is absolutely convergent when 
$|z_{1}|>|z_{2}|>0$.  The KZ equations have regular singular points. 
Using the theory of such
equations (cf. \cite{K}), (\ref{cor-fn}),
as a solution of the KZ equations in the
region given by this inequality, can be
analytically extended to a solution in the region given by the
inequality $|z_{2}|>|z_{1}-z_{2}|>0$, and the extension
can be expanded as a suitably
truncated series in rational
powers of $z_{2}$, $z_{1}-z_{2}$, $\log z_{2}$ and $\log
(z_{1}-z_{2})$. Since by the definition of intertwining operator
the original solution does not contain 
$\log z_{1}$ or
$\log z_{2}$, the extension cannot contain $\log z_{2}$ or $\log
(z_{1}-z_{2})$, and so the extension must be of the form (\ref{phyper}).
We take $N$ to be an integer such that (\ref{si})
holds for the $w_{(i)}$ and $w'_{(0)}$ above and for the 
$s_{i}$, $i=1, \dots, j$, appearing in the extension (\ref{phyper}).
Then $N$ depends only on ${\cal Y}_{1}$ and ${\cal Y}_{2}$.
(This argument is the same as that in \cite{TK} and in the proof of Theorem
3.5 in \cite{H2}.)  

For
general elements $w'_{(0)}\in L(k, \lambda_{0})'$ and $w_{(l)}\in L(k,
\lambda_{l})$, with $l=1, 2, 3$, we note that since the universal enveloping 
algebra $U(\hat{\frak g}_{-})$ is in fact generated by the
$g(-1)$ for $g\in {\frak g}$, the general elements are linear combinations
of elements of the form $g_{1}(-1)\cdots g_{m}(-1)w$ where $w$ is a
lowest conformal-weight vector and $g_{i}\in {\frak g}$. 
Using (\ref{3.1}), (\ref{3.7.1}) and (\ref{3.7.2}), we see easily that
\begin{equation}\label{f-cor-f}
\langle w'_{(0)}, {\cal Y}_{1}(w_{(1)}, x_{1})
{\cal Y}_{2}(w_{(2)}, x_{2})w_{(3)}\rangle
\end{equation}
can be written as
a linear combination of expressions of the same form but
with $w'_{(0)}\in L(\lambda_{0})^{*}$, $w_{(l)}\in L(\lambda_{l})$, 
$l=1, 2, 3$, and with Laurent polynomials in $x_{1}, x_{2}$ and 
$x_{1}-x_{2}$ (suitably expanded) as coefficients.
In particular, the product 
${\cal Y}_{1}(\cdot, x_{1}){\cal Y}_{2}(\cdot, x_{2})$ of intertwining
operators is determined by the expressions (\ref{f-cor-f}) with the 
four vectors in the lowest conformal-weight spaces. Moreover,
it is easy to see that for
general $w'_{(0)}$, $w_{(l)}$, $l=1, 2, 3$, 
(\ref{cor-fn}) is absolutely convergent when 
$|z_{1}|>|z_{2}|>0$ and can be extended to  a
multivalued analytic function of the form (\ref{phyper}) 
when $|z_{2}|>|z_{1}-z_{2}|>0$; 
note that for $n>0$, $x_{1}^{-n}$ is to be written as
$x_{2}^{-n}(1+(x_{1}-x_{2})/x_{2})^{-n}$.
In addition, each application of
$g(-1)$ to $w_{(1)}$ or $w_{(2)}$ introduces at most one expression
$(x_{1}-x_{2})^{-1}$, and each application of $g(-1)$ to $w'_{(0)}$ or
$w_{(3)}$ introduces no such poles.  We make this explicit (proving the
inequality (\ref{si})) by using induction on
the sum of the weights of $w'_{(0)}$ and the $w_{(l)}$, $l=1, 2, 3$:

Assume that the convergence and extension property holds when 
$$\wt w'_{(0)}+\wt w_{(1)}+\wt w_{(2)} +\wt w_{(3)}< q.$$
Any homogeneous element of $L(k, \lambda)$
is either a lowest-weight vector or a linear combination of elements of 
the form $g(-1)\tilde{w}$ where $\tilde{w}$ has lower weight. Thus
when 
$$\wt w'_{(0)}+\wt w_{(1)}+\wt w_{(2)} +\wt w_{(3)}=q,$$
at least one of the elements 
$w'_{(0)}, w_{(1)}, w_{(2)},  w_{(3)}$ is a linear combination
of elements of the form $g(-1)\tilde{w}$.  We shall  prove only the case
$w_{(1)}=g(-1)\tilde{w}_{(1)}$ for some 
$\tilde{w}_{(1)}\in L(k, \lambda_{1})$, the other cases being similar.

By (\ref{3.1}), 
$$
{\cal Y}_{1}(g(-1)\tilde{w}_{(1)}, x_{1})
=\nord 
g(x_{1})
{\cal Y}_{1}(\tilde{w}_{(1)}, x_{1})\nord,
$$
and so by (\ref{3.7.1}),
\begin{eqnarray}\label{4.1}
\lefteqn{\langle w'_{(0)}, {\cal Y}_{1}(w_{(1)}, x_{1})
{\cal Y}_{2}(w_{(2)}, x_{2})w_{(3)}\rangle
}\nno\\
&&=-\sum_{p\in {\Bbb N}}x_{1}^{p}\langle g(p+1)w'_{(0)}, 
{\cal Y}_{1}(\tilde{w}_{(1)}, x_{1})
{\cal Y}_{2}(w_{(2)}, x_{2})w_{(3)}\rangle
\nno\\
&&\quad +\sum_{p\in {\Bbb N}}x_{1}^{-p-1}\langle w'_{(0)}, 
{\cal Y}_{1}(\tilde{w}_{(1)}, x_{1})
{\cal Y}_{2}(w_{(2)}, x_{2})g(p)w_{(3)}\rangle
\nno\\
&&\quad +(x_{1}-x_{2})^{-1} \langle w'_{(0)}, 
{\cal Y}_{1}(\tilde{w}_{(1)}, x_{1})
{\cal Y}_{2}(gw_{(2)}, x_{2})w_{(3)}\rangle.
\end{eqnarray}
Note that the sums are finite.
Since the sum of the weights of the module elements in each term of
the right-hand side of (\ref{4.1}) is less than $q$, by the induction
assumption we obtain the desired conclusion for elements the sum of whose
weights is $q$.  This proves that products of
intertwining operators for $L(k, 0)$ have the convergence and extension
property. 

Using the KZ equations (\ref{kzn}), the same method as for products of 
two intertwining operators shows that (\ref{cor-fn-n})
is absolutely convergent when $|z_{1}|>\cdots>|z_{n}|>0$.

We now prove the general case ($m\ge 1$). Let $m>1$
and $k_{i}\in {\Bbb N}$, for
$i=1, \dots, m$. Let $W$ be a lower-truncated generalized 
module for $L(k_{1}, 0)\otimes \cdots \otimes 
L(k_{m}, 0)$.
We shall actually prove that $W$ is a direct sum of irreducible
modules, using arguments of \cite{FHL} and \cite{DMZ}.
Any element $w\in W$ 
generates a weak $L(k_{i}, 0)$-module, $i=1, \dots, m$.  Recall 
that for any $i=1, \dots, m$, any 
weak $L(k_{i}, 0)$-module is completely
reducible and the irreducible weak $L(k_{i}, 0)$-modules are in fact
the $L(k_{i}, 0)$-modules listed 
in Theorem \ref{fz}.  Thus the weak $L(k_{i}, 0)$-module
generated by $w$
is a finite direct sum of irreducible $L(k_{i}, 0)$-modules.
So the tensor product of these $L(k_{i}, 0)$-modules
for $i=1, \dots, m$,
as an $L(k_{1}, 0)\otimes \cdots \otimes L(k_{m}, 0)$-module,
is a finite direct sum of tensor products 
of irreducible $L(k_{i}, 0)$-modules
for $i=1, \dots, m$.
By Proposition 4.7.2 in \cite{FHL}, these tensor products of
irreducible modules are irreducible modules for
$L(k_{1}, 0)\otimes \cdots \otimes L(k_{m}, 0)$,
and so we have a finite direct sum of irreducible
$L(k_{1}, 0)\otimes \cdots \otimes L(k_{m}, 0)$-modules, which
canonically maps onto the generalized
$L(k_{1}, 0)\otimes \cdots \otimes L(k_{m}, 0)$-module generated by
$w$.  This in turn is a finite direct sum of modules and so 
$W$ is completely reducible.
If $W$ is finitely generated, it must be 
a module for $L(k_{1}, 0)\otimes \cdots \otimes 
L(k_{m}, 0)$.
(This argument also works using the notion of module in \cite{FZ} in
place of the notion of weak module; cf. Proposition 2.7 of \cite{DMZ}.)

Recall that for a vertex operator algebra 
$V$ and $V$-modules $W_{1}, W_{2}, W_{3}$,
the fusion rule ${\cal N}_{W_{1}W_{2}}^{W_{3}}$ is the dimension of the
space ${\cal V}_{W_{1}W_{2}}^{W_{3}}$ 
of all intertwining operators of type ${W_{3}\choose W_{1}W_{2}}$ 
(see \cite{FHL}).
By Proposition 2.10 of \cite{DMZ} (see also \cite{L}), we have the following:
Let $m>0$, $k_{i}\in {\Bbb N}$,
$V=L(k_{1}, 0)\otimes \cdots \otimes L(k_{m},
0)$, and $W_{t}=L(k_{1}, \lambda_{1}^{(t)})\otimes \cdots
\otimes L(k_{m}, \lambda_{m}^{(t)})$, for
$t=1, 2, 3$, irreducible $V$-modules (recall Theorem 4.7.4 of \cite{FHL}),
and let ${\cal Y}$ be an
intertwining operator of type ${W_{3}}\choose {W_{1}W_{2}}$. For
convenience, we shall  write the fusion rules ${\cal N}^{L(k_{i},
\lambda_{i}^{(3)})}_{L(k_{i},
\lambda_{i}^{(1)})\; L(k_{i}, \lambda_{i}^{(2)})}$, $i=1, \dots, m$,
as ${\cal N}^{3; i}_{12}$. Then
there exist intertwining operators ${\cal Y}^{(l_{i})}_{i}$ 
of type ${L(k_{i}, \lambda_{i}^{(3)})} \choose
{L(k_{i}, \lambda_{i}^{(1)})\; L(k_{i}, \lambda_{i}^{(2)})}$, 
 $l_{i}=1, \dots, 
{\cal N}^{3; i}_{12}$, $i=1, \dots, m$, such that
\begin{equation}\label{3.5}
{\cal Y}=\sum_{l_{1}=1}^{{\cal N}^{3; 1}_{12}}
\cdots
\sum_{l_{m}=1}^{{\cal N}^{3; m}_{12}}
{\cal Y}^{(l_{1})}_{1}\otimes \cdots\otimes 
{\cal Y}^{(l_{m})}_{m},
\end{equation}
where both the left- and right-hand sides are understood as
linear maps {}from 
$W_{1}\otimes W_{2}$ to $W_{3}\{x\}$.

This result reduces the convergence and extension property for
products of intertwining operators for the vertex operator algebra
$L(k_{1}, 0)\otimes \cdots \otimes L(k_{m}, 0)$ to the corresponding
properties that we have proved above (as in the proof of Theorem 3.5
in \cite{H2}).  Similarly, the third conclusion of the theorem follows
immediately {}from (\ref{3.5}) and the case $m=1$ proved above.  \epfv

\begin{rema}
{\rm The proof of the convergence and extension property for
products of two intertwining operators and the proof of the convergence of
products of an arbitrary number of
intertwining operators are basically
the same as in \cite{TK}. But Tsuchiya and Kanie 
analytically extended these convergent 
products of an arbitrary number of 
intertwining operators to the whole configuration
space and had to show that there are no logarithms of the variables occurring 
in the extensions, in order to construct
the braid group representations. In our (more general) 
theory, we need only prove the extension property
of products of two intertwining operators, and  also, this extension 
property  requires only that products can be extended to the 
region  $|z_{2}|>|z_{1}-z_{2}|>0$.
The extension property
for products of an arbitrary number of
intertwining operators to the whole configuration space
then follows {}from the convergence of the 
same number of intertwining operators, the associativity of
intertwining operators (proved using the 
extension
property for products of two intertwining operators) and skew-symmetry
and other properties of intertwining operators.}
\end{rema}

We now define the following
class of vertex operator algebras:

\begin{defi}
{\rm Let $m\in {\Bbb Z}_{+}$, $k_{i}\in {\Bbb N}$.  
A vertex operator algebra $V$ is
said to be {\it in the class ${\cal L}_{k_{1}, \dots, k_{m}}$}
if $V$ has a vertex operator subalgebra (with the same Virasoro
element as $V$)
isomorphic to $L(k_{1},
0)\otimes \cdots \otimes L(k_{m}, 0)$.}
\end{defi}

For vertex operator algebras in this class,
by Theorem 3.2 in \cite{H2} and 
Theorem \ref{kim} above, the condition for the applicability 
of the tensor product theory remaining to be 
verified is
the following:

\begin{propo}\label{kim1}
Let $V$ be a vertex operator algebra in the class 
${\cal L}_{k_{1}, \dots, k_{m}}$. Then 
every finitely-generated lower-truncated generalized $V$-module 
is a $V$-module. 
\end{propo}
\pf
The space $V$ is a module for the (rational) vertex operator algebra 
$L(k_{1}, 0)\otimes \cdots \otimes L(k_{m}, 0)$ and so is a 
finite direct sum of tensor products of  irreducible $L(k_{i},
0)$-modules by Theorem 4.7.4 of \cite{FHL}.
Since irreducible $L(k_{i}, 0)$-modules are spanned by elements of
the form $g_{1}(-1)\cdots g_{m}(-1)v$ where $v$ is a
lowest conformal-weight vector and $g_{i}\in {\frak g}$, 
$V$ is spanned by tensor products of elements of this form,
where the needed lowest conformal-weight vectors range through a
finite-dimensional subspace of $V$.
Let $W$ be a lower-truncated generalized $V$-module generated by 
a single element $w$.
Then $W$ is also a lower-truncated generalized
$L(k_{1}, 0)\otimes \cdots \otimes L(k_{m}, 0)$-module, and the
$L(k_{1}, 0)\otimes \cdots \otimes L(k_{m}, 0)$-submodule generated by
$w$ is a module by Theorem \ref{kim}.
By a lemma of Dong-Mason \cite{DM} and Li
\cite{L}, $W$ is spanned by elements of the form
$u_{n}w$ where $u\in V$, $n\in {\Bbb Z}$.
Thus {}from (\ref{3.1}) for the vertex operator map associated with
$W$, we see easily that $W$ is a module.
(Cf. the proof of Proposition 3.7 in \cite{H2}.)
\epfv

Let $V$ be a vertex operator algebra in the class ${\cal L}_{k_{1}, 
\dots, k_{m}}$.  For the 
conformal equivalence class $P(z)$
of spheres with negatively oriented puncture
$\infty$ and positively oriented punctures $z$ and $0$ and with the 
trivial local coordinates vanishing at these punctures,
recall 
the bifunctor $\hboxtr_{P(z)}$ and the tensor product bifunctor
$\boxtimes_{P(z)}$
constructed in \cite{HL5}.
By Theorem \ref{kim}, Proposition
\ref{kim1}, and Theorems 3.1 and 3.2 and Corollary 3.3 in \cite{H2}, which 
in turn are proved 
using results in \cite{HL1}--\cite{HL6} and \cite{H1}, we obtain the
following results:

\begin{propo}
For any $V$-modules $W_{1}$ and $W_{2}$,
$W_{1}\hboxtr_{P(z)}W_{2}$ is a module and the
$P(z)$-product $(W_{1}\boxtimes_{P(z)} W_{2}, Y_{P(z)};
\boxtimes_{P(z)})$ is a $P(z)$-tensor product of $W_{1}$ and $W_{2}$.\epf
\end{propo}

\begin{theo}
For any $V$-modules $W_{1}$, $W_{2}$ and $W_{3}$ and any complex
numbers $z_{1}$ and $z_{2}$ satisfying
$|z_{1}|>|z_{2}|>|z_{1}-z_{2}|>0$, there exists a unique isomorphism
${\cal A}_{P(z_{1}), P(z_{2})}^{P(z_{1}-z_{2}), P(z_{2})}$ {}from
$W_1\boxtimes_{P(z_1)}(W_2 \boxtimes_{P(z_2)}W_3)$ to $(W_1
\boxtimes_{P(z_1-z_2)}W_2)\boxtimes_{P(z_2)}W_3$ such that
for any $w_{(1)}\in W_{1}$, $w_{(2)}\in W_{2}$ and $w_{(3)}\in W_{3}$,
\begin{eqnarray*}
\lefteqn{\overline{{\cal A}}_{P(z_{1}), P(z_{2})}^{P(z_{1}-z_{2}), 
P(z_{2})}(w_{(1)}\boxtimes_{P(z_{1})}(w_{(2)}
\boxtimes_{P(z_{2})}w_{(3)}))}\nno\\
&&=(w_{(1)}\boxtimes_{P(z_{1}-z_{2})}w_{(2)})
\boxtimes_{P(z_{2})}w_{(3)},
\end{eqnarray*}
where $$\overline{{\cal A}}_{P(z_{1}), P(z_{2})}^{P(z_{1}-z_{2}),
P(z_{2})}: \overline{W_1\boxtimes_{P(z_1)}(W_2 \boxtimes_{P(z_2)}W_3)}
\to \overline{(W_1 \boxtimes_{P(z_1-z_2)}W_2)\boxtimes_{P(z_2)}W_3}$$
is the canonical extension of ${\cal A}_{P(z_{1}), 
P(z_{2})}^{P(z_{1}-z_{2}), P(z_{2})}$.\epf
\end{theo}

\begin{theo}[associativity for intertwining operators]\label{assoc}

\begin{enumerate}

\item For any $V$-modules 
$W_{0}$, $W_{1}$, $W_{2}$, $W_{3}$  and $W_{4}$, 
any 
intertwining operators ${\cal Y}_{1}$ and ${\cal Y}_{2}$ of 
 types ${W_{0}}\choose {W_{1}W_{4}}$ and ${W_{4}}\choose {W_{2}W_{3}}$,
respectively,
and any choice of $\log z_{1}$ and $\log z_{2}$,
$$\langle w'_{(0)}, {\cal Y}_{1}(w_{(1)},
x_{1}){\cal Y}_{2}(w_{(2)}, x_{2})w_{(3)}\rangle\lbar_{x^{n}_{1}
=e^{n\log z_{1}},\;
x^{n}_{2}=e^{n\log z_{2}}, \; n\in {\Bbb C}}$$
is absolutely convergent when $|z_{1}|>|z_{2}|>0$ for
$w'_{(0)}\in W'_{0}$, $w_{(1)}\in W_{1}$,
$w_{(2)}\in W_{2}$ and $w_{(3)}\in W_{3}$.
For any modules 
$W_{0}$, $W_{1}$, $W_{2}$, $W_{3}$,  and $W_{5}$ and 
any intertwining operators ${\cal Y}_{3}$
and ${\cal Y}_{4}$ of types ${W_{5}}\choose {W_{1}W_{2}}$ and 
${W_{0}}\choose {W_{5}W_{3}}$, respectively, and any choice of 
$\log z_{2}$ and $\log (z_{1}-z_{2})$,
$$\langle w'_{(0)}, {\cal Y}_{4}({\cal Y}_{3}(w_{(1)},
x_{0})w_{(2)}, x_{2})w_{(3)}\rangle\lbar_{x^{n}_{0}=e^{n \log
(z_{1}-z_{2})},\;
x^{n}_{2}=e^{n\log z_{2}}, \; n\in {\Bbb C}}$$
is absolutely convergent when $|z_{2}|>|z_{1}-z_{2}|>0$
for $w'_{(0)}\in W'_{0}$, $w_{(1)}\in W_{1}$,
$w_{(2)}\in W_{2}$ and $w_{(3)}\in W_{3}$.

\item For any $V$-modules 
$W_{0}$, $W_{1}$, $W_{2}$, $W_{3}$  and $W_{4}$ and 
any 
intertwining operators ${\cal Y}_{1}$ and ${\cal Y}_{2}$ of 
types ${W_{0}}\choose {W_{1}W_{4}}$ and ${W_{4}}\choose {W_{2}W_{3}}$,
respectively,
there exist a module $W_{5}$ and intertwining operators ${\cal Y}_{3}$
and ${\cal Y}_{4}$ of types ${W_{5}}\choose {W_{1}W_{2}}$ and 
${W_{0}}\choose {W_{5}W_{3}}$, respectively, such that
for any $z_{1}, z_{2}\in {\Bbb C}$
satisfying $|z_{1}|>|z_{2}|>|z_{1}-z_{2}|>0$ and any fixed choices of
$\log z_{1}$, $\log z_{2}$ and $\log (z_{1}-z_{2})$,
\begin{eqnarray}\label{1-1}
\lefteqn{\langle w'_{(0)}, {\cal Y}_{1}(w_{(1)},
x_{1}){\cal Y}_{2}(w_{(2)}, x_{2})w_{(3)}\rangle\lbar_{x^{n}_{1}
=e^{n\log z_{1}},\;
x^{n}_{2}=e^{n\log z_{2}}, \; n\in {\Bbb C}}}
\nno\\
&&=\langle w'_{(0)}, {\cal Y}_{4}({\cal Y}_{3}(w_{(1)},
x_{0})w_{(2)}, x_{2})w_{(3)}\rangle\lbar_{x^{n}_{0}=e^{n \log
(z_{1}-z_{2})},\;
x^{n}_{2}=e^{n\log z_{2}}, \; n\in {\Bbb C}}\nno\\
&&
\end{eqnarray}
for any 
$w'_{(0)}\in W'_{0}$, $w_{(1)}\in W_{1}$,
$w_{(2)}\in W_{2}$ and $w_{(3)}\in W_{3}$.

\item For any modules 
$W_{0}$, $W_{1}$, $W_{2}$, $W_{3}$,  and $W_{5}$ and 
any intertwining operators ${\cal Y}_{3}$
and ${\cal Y}_{4}$ of types ${W_{5}}\choose {W_{1}W_{2}}$ and 
${W_{0}}\choose {W_{5}W_{3}}$, respectively, there exist a module $W_{4}$ 
and intertwining operators ${\cal Y}_{1}$ and ${\cal Y}_{2}$ of 
types ${W_{0}}\choose {W_{1}W_{4}}$ and ${W_{4}}\choose {W_{2}W_{3}}$,
respectively, such that
for any $z_{1}, z_{2}\in {\Bbb C}$
satisfying $|z_{1}|>|z_{2}|>|z_{1}-z_{2}|>0$ and any fixed choices of
$\log z_{1}$, $\log z_{2}$ and $\log (z_{1}-z_{2})$, 
(\ref{1-1}) holds for any
$w'_{(0)}\in W'_{0}$, $w_{(1)}\in W_{1}$,
$w_{(2)}\in W_{2}$ and $w_{(3)}\in W_{3}$.\epf

\end{enumerate}
\end{theo}

\begin{theo}[commutativity for intertwining operators]
For any $V$-modules 
$W_{0}$, $W_{1}$, $W_{2}$, $W_{3}$  and $W_{4}$ and 
any 
intertwining operators ${\cal Y}_{1}$ and ${\cal Y}_{2}$ of 
types ${W_{0}}\choose {W_{1}W_{4}}$ and ${W_{4}}\choose {W_{2}W_{3}}$,
respectively,
there exist a module $W_{5}$ and intertwining operators ${\cal Y}_{3}$
and ${\cal Y}_{4}$ of types ${W_{0}}\choose {W_{2}W_{5}}$ and 
${W_{5}}\choose {W_{1}W_{3}}$, respectively, such that for any 
$w'_{(0)}\in W'_{0}$, $w_{(1)}\in W_{1}$,
$w_{(2)}\in W_{2}$ and $w_{(3)}\in W_{3}$,
the multivalued analytic function 
$$\langle w'_{(0)}, {\cal Y}_{1}(w_{(1)},
x_{1}){\cal Y}_{2}(w_{(2)}, x_{2})w_{(3)}\rangle\lbar_{x_{1}
=z_{1},\;
x_{2}=z_{2}}$$
of $z_{1}$ and $z_{2}$ in the region $|z_{1}|>|z_{2}|>0$ and 
the multivalued analytic function 
$$\langle w'_{(0)}, {\cal Y}_{3}(w_{(2)},
x_{2}){\cal Y}_{4}(w_{(1)}, x_{1})w_{(3)}\rangle\lbar_{x_{1}
=z_{1},\;
x_{2}=z_{2}}$$
of $z_{1}$ and $z_{2}$ in the region $|z_{2}|>|z_{1}|>0$ are analytic
extensions of each other.\epf
\end{theo}

In \cite{H2.5}, the notion of intertwining operator algebra
was introduced (see also \cite{H4} and \cite{H6}). By Theorem 3.5
in \cite{H2.5},  we obtain:

\begin{theo}
Assume that $V$ is rational. 
Let ${\cal A}=\{a_{i}\}_{i=1}^{m}$ 
be  the set of all equivalence 
classes of
irreducible $V$-modules.  Let $W^{a_{1}}, \dots, W^{a_{m}}$ be
repesentatives of $a_{1}, \dots, a_{m}$, respectively.  Let
$W=\coprod_{i=1}^{m}W^{a_{i}}$, and let ${\cal
V}_{a_{i}a_{j}}^{a_{k}}$, for $a_{i}, a_{j}, a_{k}\in {\cal A}$, be the space
of intertwining operators of type ${W^{a_{k}}\choose
W^{a_{i}}W^{a_{j}}}$. Then $(W, {\cal A}, \{{\cal
V}_{a_{i}a_{j}}^{a_{k}}\}, {\bf 1}, \omega)$ (where ${\bf 1}$ and
$\omega$ are the vacuum and Virasoro element of $V$) is an 
intertwining operator algebra.\epf
\end{theo}

Recall the sphere partial operad $K=\{K(j)\}_{j\in {\Bbb N}}$, the
vertex partial operads $\tilde{K}^{c}=\{\tilde{K}^{c}(j)\}_{j\in {\Bbb N}}$ 
of central
charge $c\in {\Bbb C}$ constructed in \cite{H3} and 
the definition of vertex tensor category in \cite{HL4} and 
\cite{HL6}. For any $c\in {\Bbb C}$ and $j\in {\Bbb N}$, $\tilde{K}^{c}(j)$
is a trivial holomorphic line bundle over $K(j)$ and we have a 
canonical holomorphic section $\psi_{j}$. Given a vertex tensor category,
we have, among other things, 
a tensor product bifunctor $\boxtimes_{\tilde{Q}}$ for each
$\tilde{Q}\in \tilde{K}^{c}(2)$. In particular,  
$\psi_{2}(P(z))\in \tilde{K}^{c}(2)$ and thus there is a 
tensor product bifunctor $\boxtimes_{\psi_{2}(P(z))}$.

\begin{theo}\label{vtc}
Let $c$ be the central charge of $V$.
Then the category of $V$-modules has a natural structure of  vertex
tensor category of central charge $c$ such that for each $z\in {\Bbb
C}^{\times}$, the tensor product bifunctor $\boxtimes_{\psi_{2}(P(z))}$
associated with $\psi_{2}(P(z))\in
\tilde{K}^{c}(2)$ is equal to $\boxtimes_{P(z)}$ constructed
in \cite{HL5}.\epf
\end{theo}

Combining Theorem \ref{vtc}  with Theorem 4.4 in \cite{HL4} (see
\cite{HL6} for the proof), we obtain:

\begin{corol}\label{kimbtc}
Let $V$ be a vertex operator algebra in the class
${\cal L}_{k_{1},\dots, k_{m}}$. 
Then the category of $V$-modules has a natural 
structure of  braided tensor category such that the tensor product
bifunctor is $\boxtimes_{P(1)}$. In particular, 
the category of $L(k_{1}, 0)\otimes \cdots \otimes 
L(k_{m}, 0)$-modules has a natural 
structure of   braided tensor category.\epf
\end{corol}

The special case $V=L(k, 0)$ states (recall the terminology defined 
at the beginning of the introduction):

\begin{theo}\label{affbtc}
For $k\in {\Bbb N}$,   the category
generated by the standard $\hat{\frak g}$-modules of level $k$ has a natural 
structure of  braided tensor category such that the tensor product
bifunctor is $\boxtimes_{P(1)}$.\epf
\end{theo}

{\small \sc Department of Mathematics, Rutgers University,
Piscataway, NJ 08854}

{\em E-mail address}: yzhuang@math.rutgers.edu

\vspace{1em}

{\small \sc Department of Mathematics, Rutgers University,
Piscataway, NJ 08854}

{\em E-mail address}: lepowsky@math.rutgers.edu

\end{document}